\documentclass[12pt]{amsart}
\usepackage{amscd}
\input prepictex   \input pictex    \input postpictex

\DeclareMathOperator\rep{rep}
\DeclareMathOperator\ind{ind}

\DeclareMathOperator\rad{rad}
\DeclareMathOperator\soc{soc}

\DeclareMathOperator\Hom{Hom}
\DeclareMathOperator\End{End}
\DeclareMathOperator\im{Im}
\DeclareMathOperator\incl{incl}

\renewcommand\mod{{\rm mod}}
\newcommand\sub{\subseteq}
 \def \lto#1{\;\mathop{\longrightarrow}\limits^{#1}\;} 
 \def \sto#1{\;\mathop{\to}\limits^{#1}\;}
\newtheoremstyle{mytheorems}{9pt}{6pt}{\itshape}{0pt}{\sc}{.}{ }{}
\newtheoremstyle{myremarks}{6pt}{3pt}{\normalfont}{0pt}{\it}{.}{ }{}

\theoremstyle{mytheorems}
\newtheorem{theorem}{Theorem}[section]
\newtheorem{lemma}[theorem]{Lemma}
\newtheorem{corollary}[theorem]{Corollary}
\newtheorem{proposition}[theorem]{Proposition}

\theoremstyle{myremarks}

\newtheorem{example}[theorem]{Example}

\newtheorem*{definition}{Definition}

\newcommand\comment[1]{}

\begin{document}
\phantom{\footnotesize [ps-bounded2b, February 11, 2007]}

\vspace{2cm}
\centerline{\Large Abelian groups with a $p^2$-bounded subgroup,}

\smallskip\centerline{\Large revisited}
\bigskip
\centerline{By}
\medskip
\centerline{\sc Carla Petroro and Markus Schmidmeier}

\bigskip\medskip

\centerline{\parbox{10cm}{\footnotesize{\it Abstract.}
Let $\Lambda$ be a commutative local uniserial ring of length $n$, 
$p$ a generator of the maximal ideal, and $k$ the radical factor field.
The pairs $(B,A)$ where $B$ is a finitely generated $\Lambda$-module
and $A\sub B$ a submodule of $B$ such that $p^mA=0$ form the objects in the
category $\mathcal S_m(\Lambda)$.
We show that in case  $m=2$ the categories $\mathcal S_m(\Lambda)$ are 
in fact quite similar to each other:
If also $\Delta$ is 
a commutative local uniserial ring of length $n$ and with radical factor
field $k$, then the categories 
$\mathcal S_2(\Lambda)/\mathcal N_\Lambda$ and $\mathcal S_2(\Delta)/\mathcal
N_\Delta$ are equivalent for certain nilpotent categorical ideals 
$\mathcal N_\Lambda$ and $\mathcal N_\Delta$.
As an application, we recover the known classification of all pairs $(B,A)$
where $B$ is a finitely generated abelian group and $A\sub B$ a subgroup
of $B$ which is $p^2$-bounded for a given prime number $p$.     }}

\renewcommand{\thefootnote}{}
\footnotetext{{\it 2000 Mathematics Subject Classification : 
20E07 (primary), % Subgroup theorems
16G20 (secondary) % Representations of quivers
}}
\footnotetext{{\it Keywords : subgroups of abelian groups, uniserial rings, 
    poset representations}
  }

\section{History and Introduction}\label{intro}
Let $\Lambda$ be a commutative local uniserial ring of length $n$ with
radical generator $p$ and radical factor field $k=\Lambda/p$.
We consider pairs 
$(B; A)$ where $B$ is a finitely generated $\Lambda$-module and $A$
a submodule of $B$.
Such pairs form the objects in  the category $\mathcal S(\Lambda)$; 
a morphism from $(B;A)$ to $(D;C)$ is  given by a map $f:B\to D$
which satisfies $f(A)\subset C$.
We are particularly interested in the full subcategories $\mathcal S_m(\Lambda)\subset
\mathcal S(\Lambda)$ (for $m\leq n$ a natural number)
which consist of those pairs $(B;A)$ that satisfy $p^mA=0$.
For example if $\Lambda=\mathbb Z/p^n$ then we are dealing with pairs $(B;A)$
where $B$ is a finite abelian $p^n$-bounded group and $A\sub B$ a subgroup 
satisfying $p^mA=0$. 

\smallskip
Each category  $\mathcal S_m(\Lambda)$ has the 
Krull-Remak-Schmidt property, so every object has a unique direct sum decomposition into
indecomposable ones.  
Examples for indecomposable objects are {\it pickets} which are pairs $(B;A)$
where the $\Lambda$-module $B$ itself is indecomposable, hence cyclic. 
Since $\Lambda$ is uniserial, each picket $(B;A)$ is determined uniquely by the
lengths $\ell$ and $m$ of the $\Lambda$-modules $B$ and $A$; we write $P_m^\ell=(B;A)$.

\smallskip
Clearly, the complexity of categories of type $\mathcal S_m(\Lambda)$ 
increases with $m$.  The categories $\mathcal S_0(\Lambda)$ and $\mod \Lambda$
are equivalent so the only indecomposable objects in $\mathcal S_0(\Lambda)$
are the pickets of type $P_0^\ell$. 
In the category $\mathcal S_1(\Lambda)$ (which we consider briefly in Section 3),
every indecomposable object is a picket
of type $P^\ell_0$ or $P^\ell_1$.
The category $\mathcal S_2(\Lambda)$ 
contains additional indecomposables which are not pickets;
it turns out that an invariant which has been introduced by Pr\"ufer
\cite[\S7]{p} in 1923 provides an efficient classification.

\begin{definition}
Let $B$ be a $\Lambda$-module and $a\in B$ a nonzero element. 
The {\it height exponent} of $a$ is 
$$h_B(a)\;=\;h(a)\;=\;\max\{n\in\mathbb N_0\;:\; a=p^nb \;\text{for some}\;b\in B\};$$
the {\it height sequence} $H_B(a)=(h(a),h(pa),\ldots,h(p^\ell a))$ consists of the 
height exponents of the nonzero $p$-power multiples of $a$.
\end{definition}

\begin{example}\label{exampleheight}
For pickets, the height sequence consists of consecutive numbers:
In a picket $(B;A)=P^\ell_m$ where $m>0$, 
any generator for $A$ has the height sequence
$(\ell-m,\ell-m+1,\ldots,\ell-1)$. 
\end{example}

\begin{example}\label{exampleQ}
If $m\geq 2$ then there are height sequences which cannot be realized by pickets.
The sequence  $(s-1,t-1)$ where $s<t-1$ is realized by the pair
$Q^t_s=(B;A)$ where 
$B=\Lambda/(p^t)\oplus \Lambda/(p^s)$ and $A=(p^{t-2},p^{s-1})\Lambda$.
\end{example}

\smallskip
All pickets and all indecomposables of type $Q_s^t$ have the property that the subgroup
is either zero or cyclic.  This is always the case for indecomposable objects in
$\mathcal S_2(\Lambda)$ according to \cite[Theorem 4]{hrw}:

\begin{theorem}
Each pair $(B,A)\in\mathcal S_2(\mathbb Z/p^n)$ is a direct sum of indecomposable
pairs; if $(B,A)$ is an indecomposable pair then $A$ is either zero or cyclic.
\end{theorem}

A description of the indecomposable objects in terms of standard forms of matrices
is given in \cite[Theorem~7.5]{bhw}.
It turns out that whenever the pair $(B;A)$ is indecomposable with 
$A$ nonzero, then the height sequence $H_B(a)$ of a subgroup
generator $a$ uniquely determines the isomorphism type of the given pair.
Since there are $\frac12(n^2+n)$ height sequences of length at most 2 with values
at most $n$, and since there are 
$n$ isomorphism types of indecomposable pairs $(B;A)$ where $A=0$, 
we deduce that there are in total $n+\frac12(n^2+n)=\frac12(n^2+3n)$ 
indecomposable objects in $\mathcal S_2(\Lambda)$, up to isomorphism.

\medskip
In this manuscript we recover the list of indecomposable objects in 
$\mathcal S_2(\Lambda)$ using poset representations, we demonstrate
that the list does not only not depend on the choice of the base ring $\Lambda$, 
but that in fact all the categories of type $\mathcal S_2(\Lambda)$ are related:

\smallskip
Let $\Delta$ be a second commutative local uniserial ring such that $\Lambda$ 
and $\Delta$ have the same length $n$ and isomorphic radical factor fields.
Clearly, the categories $\mathcal S_2(\Lambda)$ and $\mathcal S_2(\Delta)$ cannot be
equivalent unless the rings $\Lambda$ and $\Delta$ are isomorphic. 

\smallskip
We define categorical ideals $\mathcal N_\Lambda\subset
\mathcal S_2(\Lambda)$ and $\mathcal N_\Delta\subset \mathcal S_2(\Delta)$ which are
``large enough'' to make the factor categories equivalent,
$$\mathcal S_2(\Lambda)/\mathcal N_\Lambda\;\simeq\;
  \mathcal S_2(\Delta)/\mathcal N_\Delta,$$
and ``small enough'' so that the categories $\mathcal S_2(\Lambda)/\mathcal N_\Lambda$
and $\mathcal S_2(\Lambda)$ have the same indecomposable
objects in the sense that no nonzero object in $\mathcal S_2(\Lambda)$
is isomorphic to zero when considered as object in 
$\mathcal S_2(\Lambda)/\mathcal N_\Lambda$.

\smallskip
We would like to emphasize that only basic methods from linear algebra 
are needed to establish the well-known 
list of the indecomposable representations
of this poset, and hence to obtain the list of the indecomposable objects in 
$\mathcal S_2(\Lambda)$.

\smallskip
The results in this paper have been presented at the Miami meeting of the 
AMS in April 2006.  
They are adapted from the M.S.\ thesis of the first author written under
the supervision of the second.

\section{Poset Representations}

\label{fromsub}
We introduce the poset $\mathcal P_n$, 
define a functor $F:\mathcal S_2(\Lambda)\to \rep_k\mathcal P_n$
into the category of $k$-linear representations of $\mathcal P_n$,
and show that the representations of type $F(P_m^\ell)$ and $F(Q_s^t)$
form a full list of the indecomposable representations in the image of $F$.

\smallskip Let $\mathcal P_n$ be the following poset:
$$\mathcal P_n: \qquad 
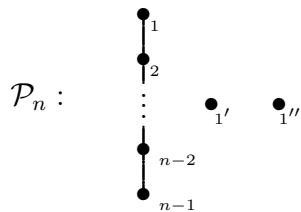
\beginpicture 
\setcoordinatesystem units <.3cm,.3cm>
\multiput {$\bullet$} at 0 4  0 2  0 -2  0 -4  3 0  6 0 /
\put{$\scriptscriptstyle 1$} at .5 3.5
\put{$\scriptscriptstyle 2$} at .5 1.5
\put{$\scriptscriptstyle n-2$} at 1.5 -2.5
\put{$\scriptscriptstyle n-1$} at 1.5 -4.5
\put{$\scriptscriptstyle 1'$} at 3.5 -.5 
\put{$\scriptscriptstyle 1''$} at 6.5 -.5
\plot 0 4  0 1  /
\plot 0 -4  0 -1 /
\put{$\vdots$} at 0 0.3 
\endpicture$$
Recall that a representation $(V^0;(V^i)_{i\in\mathcal P_n})$ for $\mathcal P_n$ is
a $k$-vector space $V^0$, the {\it total space} of $V$, together with 
subspaces $V^i\subset V^0$ for $i\in \mathcal P_n,$ such that $V^i\subset V^j$ 
holds whenever $i<j$ in $\mathcal P_n$.
For short we write $V'$ for $V^{1'}$ and $V''$ for $V^{1''}$. 

\smallskip
The category $\rep_k \mathcal P_n$ is a Krull-Remak-Schmidt category, so every
representation has a unique direct sum decomposition into indecomposable representations.
The indecomposable objects in $\rep_k\mathcal P_n$ are in $1-1$-correspondence 
to those indecomposable representations of the Dynkin diagram $\mathbb D_{n+2}$ which
have support in the central point, thus there are 
$$\#\{\ind\rep\mathbb D_{n+2}\}-\#\{\ind\rep\mathbb A_{n-1}\}-2\#\{\ind\rep\mathbb A_1\}
\quad=$$
$$\qquad\qquad=\quad(n^2+3n+2)-\frac12(n^2-n)-2\quad=\quad\frac12n^2+\frac72n$$
indecomposables. They are as follows:
The total space $V^0$ of an indecomposable representation $V$
either has dimension 1 or dimension 2.  If $\dim V^0=1$ then $V$ is isomorphic to
one of the representations $V_{\ell,\ell',\ell''}$, where 
$0\leq \ell<n$ and $0\leq \ell',\ell''\leq 1$, defined as follows. 
$$V_{\ell,\ell',\ell''}^i= 
      \left\{\!\!\begin{array}{ll}k & \text{if}\; i\leq \ell\\ 0 & \text{if}\; i>\ell \end{array}
      \right.
\;
  V_{\ell,\ell',\ell''}'= 
      \left\{\!\! \begin{array}{ll} k & \text{if}\; \ell'=1\\ 0 & \text{if}\; \ell'=0\end{array}
      \right.
\;
  V_{\ell,\ell',\ell''}''=
      \left\{\!\! \begin{array}{ll}  k & \text{if}\; \ell''=1 \\ 0 & \text{if}\;\ell''=0\end{array}
      \right. 
$$
If $\dim V^0=2$ then $V$ is isomorphic to one of the representations $W_{s,t}$
where $1\leq s<t<n$:
$$W_{s,t}^i=\left\{\begin{array}{ll} k\oplus k & \text{if}\quad i\leq s\\
                                    \Delta & \text{if}\quad s<i\leq t\\
                                    0 & \text{if}\quad i>t  \end{array}
            \right.
\qquad
  W_{s,t}'=k\oplus 0
\qquad
  W_{s,t}''=0\oplus k
$$
where $\Delta=k(1,1)\subset k\oplus k$ is the diagonal.

\bigskip 
Given a pair $(B;A)\in\mathcal S_2(\Lambda)$, we obtain a representation $V$ of $\mathcal P_n$
as follows.  Consider the filtration for $B$
given by the subspace $\rad A$:
$$\begin{array}{llll} L_0 & = A^- & =\rad A \\
                L_1 & = A^+ & =p^{-1}\rad A & =A+\soc B \\
                L_2 & & =p^{-2}\rad A \\
                & & \vdots \\
                L_{n-1} & & = p^{1-n}\rad A \\
                L_n & & = p^{-n}\rad A & = B\end{array}$$
Here, as usual in this manuscript, we write $p$ for the endomorphism of $B$
given by multiplication by $p$.  Thus, for a submodule $U\subset B$ we denote by
$pU$ and $p^{-1}U$ the image and the inverse image of $U$ under this map. 
Note that subsequent quotients of the filtration for $B$ are vector spaces; 
in particular, $V^0=A^+/A^-$ will be the total space.
For $\ell>0$, the multiplication by $p^\ell$
defines  maps $p^\ell:L_{\ell+1} \to A^+$, $p^\ell:L_\ell\to A^-$
which have image $\rad^\ell B\cap A^+$ and $\rad^\ell B\cap A^-$, respectively, 
and which give rise to  isomorphisms
$$\frac{L_{\ell+1}}{L_\ell}\lto\cong 
  \frac{\rad^\ell B\cap A^+}{\rad^\ell B\cap A^-}
  \lto\cong \frac{(\rad^\ell B\cap A^+)+A^-}{A^-}$$
into the submodule $V^\ell =((\rad^\ell B\cap A^+)+A^-)/A^-$ of $V^0$. 
We also set $V'=\soc B/ A^-$ and $V''=A/ A^-$.

\smallskip
Note that all the spaces $V^j$ have the form 
$\widetilde C\;=\;((C\cap A^+)+A^-)/A^-\;\sub \;A^+/A^-$
for a suitable submodule $C$ of $B$.  
% Namely, $V^0=\widetilde B$, $V^\ell=\widetilde{\rad^\ell B}$ for $1\leq \ell<n$,
% $V'=\widetilde{\soc B}$ and $V''=\widetilde A$. 
We collect some properties of this construction $C\mapsto \widetilde{C}$.

\newcommand\C{C^{\phantom :}}

\begin{lemma} \label{tilde} Let $(B;A)$ be a pair in $\mathcal S_2(\Lambda)$ and 
$C$, $C'$ be submodules of $B$. 
\begin{enumerate} 
\item If $C\sub C'$ then $\widetilde\C\sub \widetilde{C'}$.
\item Always, $\widetilde{C\cap C'}\sub \widetilde\C \cap \widetilde{C'}$ holds.
  If $A^-\sub C$ or $A^-\sub C'$ then equality holds.
\item Always, $\widetilde{\C}+\widetilde{C'}\sub\widetilde{C+C'}$ holds.
  We have equality if $C\sub A^+$ or $C'\sub A^+$.
\end{enumerate}
\end{lemma}

\begin{proof}
(1) If $C\sub C'$ then $(C\cap A^+)+A^-\sub (C'\cap A^+)+A^-$ holds and the assertion
follows.

(2) The inclusion $(C\cap C'\cap A^+)+A^-\sub [(C\cap A^+)+A^-]\cap[(C'\cap A^+)+A^-]$
holds always.  If $A^-\sub C$ is given, then the right hand side in the inclusion
simplifies to $(C\cap A^+)\cap[(C'\cap A^+)+A^-]$ and by the modular law 
to $(C\cap C'\cap A^+)+A^-$.

(3) Similarly, $(C\cap A^+)+(C'\cap A^+)+A^-\sub [(C+C')\cap A^+]+A^-$ holds always.
If also $C\sub A^+$ holds, then the left hand side simplifies to
$C+(C'\cap A^+)+A^-$ and by the modular law to $[(C+C')\cap A^+]+A^-$.
\end{proof}

As a consequence we obtain:

\begin{proposition} \label{prop12}

$(1)$ The assignment which maps an object $(B,A)$ in $\mathcal S_2(\Lambda)$
    to the representation $V=(V^0,(V^j))$ of $\mathcal P_n$ given by 
    $$V^0=\widetilde B;\;
    V^\ell =\widetilde{\rad^\ell B},
    \;\text{for}\; 1\leq \ell\leq n-1,\quad
    V'=\widetilde{\soc B},\quad V''=\widetilde A,$$
    defines an additive functor $F:\mathcal S_2(\Lambda)\to \rep_k\mathcal P_n$. 

\smallskip $(2)$
    Each representation $V=F(B;A)$ satisfies the following conditions:
    $$V^{n-1}\subset V'\qquad\text{and}\qquad \;V'+V''=V^0.$$
\end{proposition}

\begin{proof}
It follows from Lemma~\ref{tilde} that $V$ is a representation of $\mathcal P_n$
satisfying (2).  If $f:(B;A)\to (D;C)$ is a morphism in $\mathcal S_2(\Lambda)$,
then $f$ maps the submodules of $B$, $A^-=\rad A$, $A$, $\rad^\ell B$ (where
$0\leq\ell\leq n-1$), $\soc B$, and $A^+$
into the corresponding submodules of $D$ and hence gives rise to a map $F(f)$ between 
the representations $F(B;A)$ and $F(D;C)$. 
\end{proof}

\medskip
\begin{definition} We denote by $\rep_k' \mathcal P_n$ the full subcategory of 
$\rep_k \mathcal P_n$
consisting of those representations which satisfy the condition (2).
\end{definition}

Among the $\frac12 n^2+\frac 72n$ indecomposable representations for $\mathcal P_n$,
the two representations $V_{n-1,0,0}$ and $V_{n-1,0,1}$, and the $n-1$ representations
$W_{s,n-1}$ where $0\leq s\leq n-2$, do not have the property that $V^{n-1}\subset V'$.
The condition that $V'+V''=V^0$ excludes the $n$ representations $V_{\ell,0,0}$ where
$0\leq \ell\leq n-1$.
It turns out that all the remaining $\frac 12n^2+\frac 32n$ indecomposable representations
in $\rep_k'\mathcal P_n$ are in bijection with the indecomposable pairs in 
$\mathcal S_2(\Lambda)$ of type $P_m^\ell$ and $Q_s^t$:

\begin{proposition} \label{dense}
The functor $F$ gives rise to the following correspondence between 
pairs in $\mathcal S_2(\Lambda)$ and representations in $\rep_k\mathcal P_n$.
$$\begin{array}{lll}
  F(P_0^\ell) &= V_{\ell-1,1,0}\qquad & (1\leq \ell\leq n) \\
  F(P_1^\ell) &= V_{\ell-1,1,1} & (1\leq \ell\leq n)\\
  F(P_2^\ell) &= V_{\ell-2,0,1}& (2\leq \ell\leq n)\\
  F(Q_s^t) &= W_{s-1,t-2} & (0\leq s-1<t-2\leq n-2)
\end{array}$$
As a consequence, $F:\mathcal S_2(\Lambda)\to \rep_k'\mathcal P_n$ is a dense functor.
Moreover, all pairs of type $P_m^\ell$ and $Q_s^t$ are indecomposable
and pairwise nonisomorphic.
\qed
\end{proposition}

We will see in Section~\ref{construct} that $F$ is full.  Hence the pairs $P_m^\ell$ and
$Q_s^t$ form a full list of the indecomposable objects in $\mathcal S_2(\Lambda)$.

\section{Picket Decomposition} \label{picket}

We show in this section that every 
object in $\mathcal S_1(\Lambda)$ is direct sum of pickets
of type $P^\ell_0$ or $P^\ell_1$ where $1\leq\ell\leq n$. 
We deduce in Theorem~\ref{picketsum} that 
every pair $(B;A)\in\mathcal S_2(\Lambda)$ with the additional property
that $\soc B\subset A$ is a direct sum of pickets of type 
$P^\ell_1$ or $P^\ell_2$.
 
\smallskip
The ring $\Lambda$ is selfinjective of Loewy length $n$.
It turns out that the pair $P_0^n=(\Lambda;0)$ 
is a relatively injective
indecomposable object in $\mathcal S(\Lambda)$ with source map
the inclusion $u_n: P_0^n\to P_1^n$.  We give a direct proof of this result
which follows also from \cite[Proposition~1.4]{rs}.

\begin{lemma}\label{3.1}
Let $(B,A)$ be an object in $\mathcal S(\Lambda)$ and
$f:P^n_0\to (B,A)$ a morphism.  Either $f$ is a split monomorphism or else
$f$ factors over the inclusion $u_n:P_0^n\to P_1^n$. 
Thus, $P_0^n$ is relatively injective in the sense that
every monomorphism from $P_0^n$ with cokernel in $\mathcal S(\Lambda)$
is a split monomorphism.
\end{lemma}

\begin{proof}
Suppose a morphism $P_0^n\to (B,A)$ is given by a map $f:\Lambda\to B$. 
Clearly, if $f$ is not a monomorphism
or if $\im f\cap A\neq 0$, then the map $P_0^n\to (B;A)$ 
factors over $u_n$.  It remains to deal with the case that $f$ is a 
monomorphism such that $E\cap A=0$ where $E=\im f$.  Let $e_A:A\to E(A)$ and 
$e_B:B\to E(B)$ be injective envelopes, then there is a map $h$ which makes
the following diagram commutative.
$$\begin{CD}
 E\oplus A @> \incl >> B\\ @V 1\oplus e_A VV @VV e_B V \\ E\oplus E(A)@>h>> E(B)
\end{CD}$$
Since $E\oplus A$ is large in $E\oplus  E(A)$ and since the composition
$h\circ(1\oplus e_A)$ is a monomorphism, also $h$ is a monomorphism. Then
$h$ is a split monomorphism, so we can write $E(B)=E\oplus E(A)\oplus E''$
and define $B'=\{b\in B:e_B(b)\in E(A)\oplus E''\}$.  Since $E\subset B$, the
assertions $B=E\oplus B'$ and 
$A\sub B'$ hold. The decomposition $(B,A)=(E,0)\oplus (B',A)$ 
demonstrates that $P_0^n\to (B;A)$ is a split monomorphism.

Suppose a monomorphism $P_0^n\to (B;A)$ is such that the cokernel $(D;C)$
is an object in $\mathcal S(\Lambda)$.  Then we have a commutative diagram with
exact rows in which the vertical maps are monomorphisms.
$$\begin{CD}
 0 @>>> 0       @>>>  A @>>> C @>>> 0 \\ @. @VVV @VV\incl V @VV\incl V @. \\
 0 @>>> \Lambda @>f>> B @>>> D @>>> 0
\end{CD}$$
By the snake lemma, the composition $\Lambda\sto fB\to B/A$ is a monomorphism,
hence $\im f\cap A=0$ and we have just seen that this implies that $P_0^n\to (B;A)$ 
is a split monomorphism.
\end{proof}

\begin{lemma}\label{3.2}
The object $P_1^n$ is injective in $\mathcal S_1(\Lambda)$.
\end{lemma}

\begin{proof}
Let $P_1^n\to (B;A)$ be a monomorphism in $\mathcal S_1(\Lambda)$.
If the image is $(E,\soc E)$ then the injective $\Lambda$-module $E$ has a
complement in $B$, say  $B=E\oplus B'$.  Since $A$ is semisimple and 
$\soc E\sub A$ we obtain
$A\cap (E\oplus B')=A\cap (\soc E\oplus \soc B') = 
\soc E\oplus (A\cap \soc B')=\soc E\oplus (A\cap B')$ by the modular law 
and hence the pair $(B;A)$
decomposes as $(E;\soc E)\oplus (B';A\cap B')$.  
\end{proof}

\begin{proposition}\label{decomp}
Every object in $\mathcal S_1(\Lambda)$ is a direct sum of 
objects of type $P_0^\ell$ or $P_1^\ell$ where $1\leq \ell \leq n$.
\end{proposition}

\begin{proof}
Assume that the pair $(B;A)\in \mathcal S_1(\Lambda)$ is nonzero.
We show that $(B;A)$ has a summand of type $P_0^\ell$ or $P_1^\ell$ where
$\ell$ is the Loewy length of $B$.  
Indeed $B$, having Loewy length $\ell$, has an indecomposable summand, say
$E$, which is a projective-injective $\Lambda/(p^\ell)$-module.
According to Lemma~\ref{3.1}, either $P_0^\ell=(E;0)$ splits off 
as a direct summand of $(B;A)$, or else $P_1^\ell=(E;\soc E)$ can be embedded into
$(B;A)$.  In the second case, $P_1^\ell$ splits off as a direct summand since it is an
injective module (Lemma~\ref{3.2}). 
\end{proof}

\begin{theorem} \label{picketsum}
Every pair $(B;A)$ in $\mathcal S_2(\Lambda)$ with the extra property that
$\soc B\subset A$, is isomorphic to a direct sum of pickets $P_m^\ell$ where $m=1$ or $2$ 
and $m\leq \ell\leq n$. 
\end{theorem}

\begin{proof}
For each pair $(B;A)\in\mathcal S_2(\Lambda)$ which satisfies $\soc B\subset A$ consider
the picket decomposition of the corresponding pair $(B;\rad A)$ in $\mathcal S_1(\Lambda)$ 
given by Proposition~\ref{decomp},
$$(B;\rad A)\;=\;\bigoplus_{i=1}^s(B_i;\rad A\cap B_i).$$
We show that this yields the picket decomposition
for $(B;A)$:
$$ (B;A)\;=\;\bigoplus_{i=1}^s(B_i;A\cap B_i).$$
Consider $p$ as an endomorphism of $B$ and write $pU$ and
 $p^{-1}U$ for the image and the inverse image of the submodule $U\subset B$.
\begin{eqnarray*}
A &\stackrel{(1)}=& A + \soc B \\
  &=& p^{-1}pA\\
  &\stackrel{\ref{decomp}}=& p^{-1}\Big( \sum_i pA\cap B_i \Big)\\
  &\stackrel{(2)}=& \sum_i p^{-1}\left(pA\cap B_i\right) \\
  &=& \sum_i \left(p^{-1}pA\cap p^{-1}B_i\right)\\
  &\stackrel{(1)}=& \sum_i \left(A\cap (B_i+\soc B)\right)\\
  &\stackrel{(1)}=& \soc B+\sum_i \left(A\cap B_i\right)\\
  &=& \sum_i \left( \soc B_i+ A\cap B_i\right)\\
  &=& \bigoplus_i A\cap B_i
\end{eqnarray*}
In the equations labelled $(1)$, equality holds since $\soc B\subset A$,
and equality in $(2)$ follows from $pA\cap B_i\subset pB$.
\end{proof}

\section{Homomorphism Categories}\label{construct}

We show in Theorem~\ref{full} 
that the functor $F:\mathcal S_2(\Lambda)\to \rep_k\mathcal P_n$ is full.

\begin{definition}
Let $\mathcal C$ be any category.  
The {\it homomorphism category} $\mathcal H(\mathcal C)$ 
has as objects all triples
$(B;A;f)$ or $(B\leftarrow^fA)$ where $f:A\to B$ is a morphism 
in $\mathcal C$. A morphism in 
$\mathcal H(\mathcal C)$ from $(B;A;f)$ to $(D;C;e)$ 
is a pair $(h:B\to D,g:A\to C)$ of morphisms
in $\mathcal C$ which satisfies $eg=hf$. 
If $\mathcal C$ is an abelian category, then so is $\mathcal H(\mathcal C)$;
in this case kernels and cokernels, and hence pull-backs and push-outs,
are computed componentwise.
For short we write $\mathcal H(\Lambda)$ for $\mathcal H({\rm mod}\Lambda)$,
thus $\mathcal S(\Lambda)\subset \mathcal H(\Lambda)$ is an embedding 
of a full subcategory. 
\end{definition}

\begin{proposition} \label{Apm}
Given a pair $(B;A^+)\in \mathcal S_2(\Lambda)$ where $\soc B\subset A^+$,
and a subspace $U$ of $A^+/A^-$ where $A^-=\rad A^+$, 
then there is a unique submodule $A$ of $B$ such that the diagram
$$(*)\quad\begin{CD} 
0 @>>> (B;A^-) @>>> (B;A) @>>> (0; U) @>>> 0 \\
@. @| @VVV @VVV @. \\
0 @>>> (B;A^-) @>>> (B;A^+) @>>> (0; A^+/A^-) @>>> 0
\end{CD}$$
is a pull-back diagram in the category $\mathcal H(\Lambda)$. 
Conversely, every pair $(B;A)$ arises in this way. 
\end{proposition}

\begin{proof}
The bottom sequence in $(*)$ is a short exact sequence in the category $\mathcal H(\Lambda)$, 
and the pull-back
along the inclusion $U\to A^+/A^-$ yields an object $X\in\mathcal H(\Lambda)$ 
and two monomorphisms $(B;A^-)\to X\sto v (B;A^+)$.   Thus, 
$X$ is in $\mathcal S_2(\Lambda)$ and by identifying $X$ with $\im v$ we obtain
a uniquely determined submodule $A\subset B$ such that the above diagram $(*)$
is commutative and has exact rows.  

In order to realize a pair $(B;A)\in \mathcal S_2(\Lambda)$ as a 
pull-back, take $A^+=A+\soc B$, $A^-=\rad A$ and for $U$ the subspace
$A/A^-$ of $A^+/A^-$.
Then the diagram $(*)$ is commutative with exact rows and hence is 
a pull-back diagram since the vertical map on the left 
is an isomorphism.
\end{proof}

\smallskip
Note that the last term $A^+/A^-$ in the bottom sequence in $(*)$ is the total space 
of the poset representation $F(B;A^+)$.  Our next result shows that
any morphism $F(B;A^+)\to F(D;C^+)$ between poset representations can be lifted to 
a map $(B;A^+)\to (D;C^+)$:

\begin{proposition} \label{lifting}
Suppose that the pairs $(B;A)$ and $(D;C)$ in $\mathcal S_2(\Lambda)$
satisfy $\soc B\subset A$ and $\soc D\subset C$.  
Then the functor $F$ induces an isomorphism
$$\frac{\Hom_{\mathcal S_2(\Lambda)}((B;A),(D;C))}{\mathcal N((B;A),(D;C))} 
  \quad\lto{\displaystyle \cong} \quad
  \Hom_{\mathcal \rep_k\mathcal P_n} (F(B;A),F(D;C))$$
where $\mathcal N((B;A),(D;C))$ consists of 
all maps $f:(B;A)\to (D;C)$ such that $f(A)\subset \rad C$. 
\end{proposition}

\begin{proof}
According to Theorem~\ref{picketsum}, the pairs $(B;A)$ and $(D;C)$ are 
direct sums of pickets of type $P_1^s$ or $P_2^s$.  Since $F$ is an additive
functor, we may assume that both $(B;A)$ and $(D;C)$ are in fact such pickets.

The kernel of the map $F_{(B;A),(D;C)}$ is $\mathcal N((B;A),(D;C) )$,
so it remains to show that any nonzero map $h:F(B;A)\to F(D;C)$
lifts to a map $f:(B;A)\to (D;C)$ which makes the following diagram commutative.
$$\CD 0 @>>> (B;A^-) @>>> (B;A) @>>> (0;A/A^-) @>>> 0\\
 @.  @V f^- VV  @V f VV  @VV(0;h^0) V \\
 0 @>>> (D;C^-) @>>> (D;C) @>>> (0;C/C^-) @>>> 0
  \endCD 
$$
Since $\soc B\subset A$, the space $F(B;A)''=A/A^-$ is nonzero and hence 
$F(B;A)$ is a poset representation of type $V_{\ell,\ell',1}$. 
Also $F(D;C)$ must be isomorphic to a representation
of type  $V_{m,m',1}$.
If the space $V_{m,m',1}'$ is zero (that is, if $m'=0$), 
then also the space $V_{\ell,\ell',1}'$ must be zero. We obtain that the length of $A$,
which is $L=2-\ell'$, is at least the length $2-m'$ of $C$.
Using the projectivity of $A$ as a $\Lambda/(p^L)$-module, we obtain a lifting 
$g:A\to C$ for $h^0$:
$$\CD A @>>> A/A^- \\
  @V g VV @VV h^0 V \\
  C @>>> C/C^- 
\endCD $$
Similarly, whenever a space $V_{m,m',1}^i$ is zero, then so is the 
space $V_{\ell,\ell',1}^i$, and hence $\ell\leq m$ holds.  
Let $D\to E$ be the inclusion into an injective envelope.
Then we can extend the composition 
$A\sto g C\sto\incl D\to E$ to a map $e:B\to E$.  Since $\ell\leq m$,
the image of $e$ is contained in $D$ so $f=e|_{B,D}$ is an extension
of $g$ which makes the following diagram commutative:
$$\CD A @>>> B \\
   @V g VV  @VV f V \\
   C @>>> D 
\endCD $$
This extension satisfies $F(f)=h$, finishing the proof.
\end{proof}

\medskip

\begin{theorem}\label{full}
The functor $F:\mathcal S_2(\Lambda)\to \rep_k \mathcal P_n$ is full.
\end{theorem}

\begin{proof}
Given two objects $(B;A)$ and $(D;C)$ in $\mathcal S_2(\Lambda)$, 
put $U=F(B;A)$ and $V=F(D;C)$ and let $h:U\to V$ be a morphism in 
the category $\rep_k\mathcal P_n$.  We show that there is a morphism
$f:(B;A)\to (D;C)$ in $\mathcal S_2(\Lambda)$ such that $F(f)=h$.

\smallskip
Note that the representations $F(B;A)$ and $F(B;A^+)$ coincide in all 
positions, only $F(B;A)''$ may be a proper subspace of $F(B;A^+)''$;
in fact, $F(B;A^+)''$ coincides with the total space $F(B;A^+)^0=F(B;A)^0$.
Thus, $h$ defines a morphism $F(B;A^+)\to F(D;C^+)$. 
By Proposition~\ref{lifting}, this morphism lifts to a map
$f^+:(B;A^+)\to (D;C^+)$ which 
makes the diagram in $\mathcal H(\Lambda)$ commutative.
$$\CD 0 @>>> (B;A^-) @>>> (B;A^+) @>>> (0;U^0) @>>> 0 \\
 @. @V f^- VV  @V f^+ VV @VV (0,h^0) V \\
0@>>> (D;C^-) @>>> (D;C^+) @>>> (0;V^0) @>>> 0
\endCD$$
Here $f^-=f^+|_{(B;A^-),(D;C^-)}$ is the restriction of $f^+$. 
The morphism $h:F(B;A)\to F(D;C)$  also yields the following commutative diagram.
$$\CD (0;U'') @> (0;h'') >> (0;V'') \\
 @VVV @VVV \\
(0;U^0) @>> (0;h^0) > (0;V^0)
\endCD $$
We obtain the desired map $f$ as a pull-back in the abelian category 
$\mathcal H(\mathcal H(\Lambda))$ which has as objects all 
homomorphisms in 
$\mathcal H(\Lambda)$.  The two diagrams above form the following diagram
in $\mathcal H(\mathcal H(\Lambda))$, which has an exact row.
$$\CD @. @. @. (0;h'') \\
@. @. @. @VVV \\
0 @>>> f^- @>>> f^+ @>>> (0;h^0) @>>> 0
\endCD $$
The pullback of this diagram,
$$(**) \hspace{2cm} \CD 0 @>>> f^- @>>> f @>>> (0;h'') @>>> 0 \\
 @. @| @VVV @VVV \\
 0 @>>> f^- @>>> f^+ @>>> (0;h^0) @>>> 0
\endCD $$
yields an object $f$ in $\mathcal H(\mathcal H(\Lambda))$.

\smallskip
Since pull-backs are computed componentwise, $f$ is in fact a 
morphism between the two pull-backs $(B;A)$ and $(D;C)$ computed in the
category $\mathcal H(\Lambda)$ (see Proposition~\ref{Apm}). 
In other words, $f^+:(B;A^+)\to (D;C^+)$ restricts to a map
$f:(B;A)\to (D;C)$.  
$$\CD (B;A^-) @>>> (B;A) @>>> (B;A^+) \\
 @V f^- VV @V f VV @VV f^+ V \\
  (D;C^-) @>>> (D;C) @>>> (D;C^+)
\endCD $$

\smallskip 
Since $\rad A= \rad A^+$ and $\rad C=\rad C^+$ we have that 
the representations $F(B;A)$ and $F(B;A^+)$, and 
also the representations $F(D;C)$ and $F(D;C^+)$, coincide in all positions
except possibly at $1''$.  Since $f$ is the restriction of $f^+$,
Proposition~\ref{lifting} implies that the linear maps
$ F(f)^j$ and $h^j$  coincide for all $j\in\mathcal P_n\backslash\{1''\}$.
For the space at $1''$ consider the submodule component of the top row in $(**)$:
$$\CD 0 @>>> A^- @>>> A @>>> U'' @>>> 0\\
@. @V f|_{A^-,C^-} VV @V f|_{A,C} VV @VV h'' V \\
0 @>>> C^- @>>> C @>>> V'' @>>> 0
\endCD $$
Since $F(B;A)''=A/A^-=U''$ and $F(D;C)''=C/C^-=V''$ we obtain that
$F(f)''=h''$, finishing the proof that $F(f)=h$.
\end{proof}

\section{Corollaries}

We have seen in Proposition~\ref{dense} and Theorem~\ref{full} that the functor 
$F:\mathcal S_2(\Lambda)\to \rep_k'\mathcal P_n$ is full and dense.  
We first consider the kernel.

\medskip
For pairs $(B;A)$ and $(D;C)$ in $\mathcal S_2(\Lambda)$ define the following subgroup
of $\Hom_{\mathcal S}((B;A),(D;C))$:
$$\mathcal N\big((B;A),(D;C)\big) 
 \quad =\quad \big\{ f:(B;A)\to (D;C) \;\big|\; f(A^+)\subset C^-\big\}.$$
This generalizes the definition given in Proposition~\ref{lifting}.
The collection $\mathcal N_\Lambda$ or $\mathcal N$ 
of all such subgroups forms a categorical ideal
in $\mathcal S_2(\Lambda)$.

\begin{lemma}
If the length of $\Lambda$ is $n>1$ then $\mathcal N$ has nilpotency index $n+1$.
\end{lemma}

Note that if $n=1$ then $\mathcal N=0$.

\begin{proof}
Given objects $(B_i;A_i)\in\mathcal S_2(\Lambda)$
for $0\leq i\leq n+1$, morphisms $f_i:(B_{i-1};A_{i-1})\to (B_i;A_i)$ in $\mathcal N$
for $1\leq i\leq n+1$, and an element $b\in B_0$. Then $p^{n-1}b\in\soc B_0$,
hence $f_1(p^{n-1}b)\in\rad A_1$.  
For each $i$ we have that if $f_i\cdots f_1(p^{n-i}b)\in \rad A_i$ then
$f_i\cdots f_1(p^{n-i-1}b)\in A_i+\soc B_i$ and hence 
$f_{i+1}\cdots f_1(p^{n-i-1}b)\in \rad A_{i+1}$. 
Thus, $f_n\cdots f_1(b)\in\rad A_n$ and hence $f_{n+1}\cdots f_1(b)=0$.
Conversely, if $n>1$ then the nilpotency index is not less than $n+1$ since the following
composition of $n$ maps is nonzero.
$$P_1^{n}\sto\incl P_2^n \sto{\cdot p} P_2^n\sto{\cdot p} \quad \cdots \quad 
  \sto{\cdot p} P_2^{n}$$
\end{proof}

\begin{corollary}
The functor $F:\mathcal S_2(\Lambda)\to \rep_k\mathcal P_n$
induces an equivalence of categories
$$\bar F:\mathcal S_2(\Lambda)/\mathcal N\to \rep_k'\mathcal P_n.$$
\end{corollary}

\begin{proof}
We show that the kernel of $F$ is $\mathcal N$.
Let $f:(B;A)\to (D;C)$ be a map in $\mathcal S_2(\Lambda)$ and denote by
$f^+$ the map $f^+=f:(B;A^+)\to (D;C^+)$.
Since $\rad A=\rad A^+$ and $\rad C=\rad C^+$ hold, 
the linear maps $F(f)^0$, $F(f^+)^0$ coincide.
By Proposition~\ref{lifting}, $F(f^+)^0=0$ if and only if 
$f^+\in\mathcal N((B;A^+),(D;C^+))$.  Then $f(A^+)\subset C^-$, but this is the
condition for $f$ to be in $\mathcal N((B;A),(D;C))$.
\end{proof}

\medskip
Since $\mathcal N$ is a nilpotent ideal, the canonical functor 
$\mathcal S_2(\Lambda)\to \mathcal S_2(\Lambda)/\mathcal N$ preserves indecomposable
objects and reflects isomorphisms, so the isomorphism classes of objects in
$\mathcal S_2(\Lambda)$ and $\mathcal S_2(\Lambda)/\mathcal N$ are in a natural
bijection.  
As a consequence, if also $\Delta$ is a commutative local uniserial ring of length
$n$ and with radical factor field $k$, then the categories $\mathcal S_2(\Lambda)$ and 
$\mathcal S_2(\Delta)$ admit a ``natural'' bijection between their objects.

\begin{corollary}
Suppose that $\Lambda$ and $\Delta$ are commutative local uniserial rings
of the same length $n$ and with radical factor fields isomorphic to $k$. 
Then the following categories are equivalent:
$$\mathcal S_2(\Lambda)/\mathcal N_\Lambda\quad\cong\quad 
  \mathcal S_2(\Delta)/\mathcal N_\Delta \qed$$
\end{corollary}

For example, the categories $\mathcal S_2(\mathbb Z/(p^n))/\mathcal N$
and $\mathcal S_2(k[T]/(T^n))/\mathcal N$ are equivalent if $k=\mathbb Z/p$. 

\medskip
{\it Question:} We have seen that $\mathcal N$ is an ideal of nilpotency
index $n+1$.  Is there a pair $\mathcal L_\Lambda$, $\mathcal L_\Delta$,
of ideals of nilpotency index $n$ such that the above
Corollary holds with $\mathcal N$ replaced by $\mathcal L$?  
For example, an ideal of nilpotency index $n$ is given by all maps which factor
through the multiplication by $p$.
Note that we cannot expect to have an ideal $\mathcal L$
of nilpotency index $r<n$:  The endomorphism ring of the pair $X=(\Lambda;0)$
has filtration 
$$0=\mathcal L^r(X,X) \subset \mathcal L^{r-1}(X,X)\subset \cdots \subset
    \mathcal L(X,X)\subset \Lambda$$
in which not all subsequent factors can be semisimple.  Hence the factor ring
$\Lambda/\mathcal L(X,X)$ is not a field.  In particular if $\Lambda=\mathbb Z/(p^n)$
and $\Delta=k[T]/(T^n)$ then the following two endomorphism rings are not isomorphic:
$\End_{\mathcal S_2(\Lambda)/\mathcal L}(\Lambda,0)\not\cong 
  \End_{\mathcal S_2(\Delta)/\mathcal L} (\Delta,0)$.

\medskip
We can also deal with the case where $\Lambda$ is a PID and $p$ a prime element. 
By primary decomposition, any pair $(B;A)$ where $B$ is a finitely generated 
$\Lambda$-module and $A$ a $p^2$-bounded submodule of $B$ has a unique direct
sum decomposition into pair of type $(F;0)$ where $F$ is a finitely generated
free $\Lambda$-module, a pair of type $(D;0)$ where $p$ acts as automorphism on $D$,
and a pair $(B;A)$ where $B$ is $p^n$-bounded for some natural number $n$.

\begin{corollary}
Let $\Lambda$ be a principal ideal domain, $p$ a prime element, $B$ a finitely 
generated $\Lambda$-module and $A$ a submodule of $B$ which is $p^2$-bounded.
Then the pair $(B,A)$ has a direct sum decomposition, unique up to isomorphy and 
reordering, into finitely many indecomposable pairs (a) of type
$(\Lambda/(q);0)$ where $q$ is a prime power relatively prime to $p$, 
or (b) of type $P^\ell_m$ or $Q_s^t$ in $\mathcal S_2(\Lambda/(p^n))$ 
for some $n\in\mathbb N$, or (c) indecomposable projective of type $(\Lambda;0)$ \qed
\end{corollary}

Returning to Pr\"ufer's height sequences, we see that an indecomposable pair
is either given by a power of a prime ideal in $\Lambda$, or else by the height
sequence of a subgroup generator.

\begin{corollary}
Let $\Lambda$ be a principal ideal domain, $p$ a prime element, $B$ a finitely 
generated $\Lambda$-module and $A$ a submodule of $B$ which is $p^2$-bounded.
An indecomposable pair $(B;A)$ is 
\begin{enumerate}\item either isomorphic to $(\Lambda/Q;0)$ for a uniquely determined
  power $Q$ of a prime ideal, or else,
\item if $A=a\Lambda$ is (nonzero) cyclic, determined uniquely, up to isomorphism, 
  by the height   sequence $H_B(a)$. 
\end{enumerate}
Conversely, every power of a prime ideal, and every height sequence of length at most $2$,
can be realized by an indecomposable pair. \qed
\end{corollary}

\medskip
{\it Address of the authors:}
\nopagebreak

\smallskip\footnotesize
  Mathematical Sciences, 
  Florida Atlantic University,
  Boca Raton, FL 33431-0991\newline
  E-mail:   {\tt petroro@math.fau.edu} (C.P.), 
  {\tt markus@math.fau.edu} (M.S.)

\end{document}